\documentstyle[twoside,11pt]{article}

\input{mssymb}
\makeatother

\newcommand{\dis}{\displaystyle}

\newcommand{\K}{{\Bbb K}}
\newcommand{\N}{{\Bbb N}}
\newcommand{\R}{{\Bbb R}}
\newcommand{\C}{{\Bbb C}}

\newcommand{\mysec}[2]{\section*{\normalsize\hfil\sc {#1}. {#2}\hfill}%
\noindent\setcounter{theo}{0}\setcounter{section}{#1}%
\typeout{#1. #2}}%

\newcommand{\proof}{\par\noindent{\em Proof. }}

\newtheorem{theo}{Theorem}[section]
\newtheorem{lemma}[theo]{Lemma}

\newtheorem{definition}[theo]{Definition}
\newenvironment{defi}{   % ist nicht kursiv
   \begin{definition}
   \begin{em}}{
   \end{em}
   \end{definition}}

\newcounter{abc}   % Counter fr statements-environment wird deklariert
\newcounter{iiiii} % Counter fr aequivalenz-environment wird deklariert

\newenvironment{aequivalenz}
{\setcounter{iiiii}{0}
\begin{list}%
{{\rm (\roman{iiiii})}}%  Falls die items nicht angegeben sind: i)u.s.w.
{\usecounter{iiiii}
\parsep=0pt plus 1pt
\topsep=1pt plus 2pt minus 1pt
\itemsep=1pt plus 2pt minus 1pt
\leftmargin=3\baselineskip
\labelsep=.6\baselineskip
\labelwidth=2.4\baselineskip
\rightmargin 0pt}%
}%               Das war das zweite Argument von "newenvironment"
{\end{list}}

\newenvironment{statements}%
{\setcounter{abc}{0}
\begin{list}%
{{\rm (\alph{abc})}}%  Falls die items nicht angegeben sind: (a) u.s.w.
{\usecounter{abc}
\parsep=0pt plus 1pt
\topsep=1pt plus 2pt minus 1pt
\itemsep=1pt plus 2pt minus 1pt
\leftmargin=3\baselineskip
\labelsep=.6\baselineskip
\labelwidth=2.4\baselineskip
\rightmargin 0pt}%
}%               Das war das zweite Argument von "newenvironment"
{\end{list}}

\newenvironment{punkte}%
{\begin{list}%
{$\bullet$}%  Falls die items nicht angegeben sind: Punkte
{\leftmargin=3\baselineskip
\labelsep=\baselineskip
\labelwidth=2.5\baselineskip
\rightmargin 0pt}%
}%               Das war das zweite Argument von "newenvironment"
{\end{list}}

\newcommand{\myforall}{\mathop{\forall}}
\newcommand{\myexists}{\mathop{\exists}}
\newcommand{\eps}{\varepsilon}
\newcommand{\begsta}{\begin{statements}}
\def\endsta{\end{statements}}
\newcommand{\begaeq}{\begin{aequivalenz}}
\def\endaeq{\end{aequivalenz}}
\newcommand{\begpu}{\begin{punkte}}
\def\endpu{\end{punkte}}

\begin{document}
\begin{center}
{\Large\bf New proofs of Rosenthal's $\ell^{1}$--theorem and the
Josefson--Nissenzweig theorem}\\[33pt]
{\sc Ehrhard Behrends}
\end{center}
\bigskip
\begin{quotation}
\small\noindent
{\sc Abstract.}		
We give elementary proofs of the theorems mentioned in
the title. Our methods rely on a simple version of Ramsey theory and
a martingale difference lemma.

They also provide quantitative results: if a Banach
space contains $\ell^{1}$  only with a bad constant then every
bounded sequence admits a subsequence which is ``nearly'' a weak
Cauchy sequence.
\end{quotation}
\bigskip
%
%
%	Section 1
%
%
\mysec{1}{Introduction}%
The aim of this paper is to give new, elementary proofs of the
following theorems:
\begin{quote}
{\bf Rosenthal's $\ell^{1}$--theorem} ([6]): Let  $(x_{n})$  be a
bounded sequence in a Banach space  $X $.  Either there is a
subsequence which is equivalent to the  $\ell^{1}$-basis or there is
a subsequence  $(x_{n_{k}})$  which is weakly Cauchy (i.e.
$(x'(x_{n_{k}}$)) converges for every  $x'\in X'$).

{\bf The Josefson--Nissenzweig theorem} ([4], [5]): For every
in\-finite--dimensional Banach space  $X$  there are 
$x'_{1},x'_{2},\ldots$  in  $X'$  such that  $\|x'_{n}\| = 1$  for
all  $n$, but   $x'_{n} \to 0$ w.r.t. the weak*--topology.
\end{quote}

They are treated in several textbooks (see e.g. [2] where the
reader also may find some historical remarks). We are going to
present new proofs which are more elementary than those in the
literature. Also,  in the case of Rosenthal's theorem, we are able to
derive a quantitative variant which reads as follows for real spaces
(cf. section 3 where also the complex case is discussed):
If $(x_{n})$  is bounded, $\eps > 0 $,  and if for every infinite 
$M \subset \N$ there are  $i_{1} < \ldots < i_{r}$ in  $M$  and 
$a_{1},\ldots,a_{r} \in \R$ with   $\sum |a_{\rho}| = 1$  and 
$\|\sum a_{\rho}x_{i_{\rho}}\| \le \eps $,  then there is a
subsequence  $(x_{n_{k}})$  such that 
$\limsup_{k}
x'(x_{n_{k}}) -
\liminf_{k} x'(x_{n_{k}}) \le 2 \eps$
for all  $x'$  with  $\|x'\| = 1 $.

We prepare the proof of this in
{\it section 2} where we treat an easy--to--formulate Ramsey--type
theorem (the proof of which is elementary but technical).  In {\it
section 4}  we show that a slightly modified assumption in Rosenthal's
theorem even gives rise to norm convergent subsequences; this has
already been published in [1], the present argument, however, is much
simpler.

Next, in {\it section 5}, we are going to prepare
the proof of the Josefson--Nissenzweig theorem (the section contains some 
facts about Banach limits and a simple martingale lemma).  {\it Section 6} 
uses the circle of ideas presented in section~3
to prove the Josefson--Nissenzweig
theorem.  This idea is also fundamental in [2], our proof however is
more direct and gives a slightly sharper result.

\mysec{2}{An elementary Ramsey theorem}%
\begin{theo} For  $r \in \N$  let  $T_{r}$  be a family of
$r$--tupels of increasing integers.  Suppose that
$$
\myforall_{
\scriptstyle M \subset \N \atop \scriptstyle M\ \rm infinite}
\myexists_{
\scriptstyle i_{1},i_{2},\ldots\in M \atop \scriptstyle 
i_{1} < i_{2} < \ldots}
\myforall_{r}\ 
(i_{1},\ldots,i_{r}) \in T_{r}.  \leqno(1)
$$
Then it follows that
$$
\myexists_{
\scriptstyle M_{0} \subset \N \atop \scriptstyle M_{0}\ \rm infinite}
\myforall_{
\scriptstyle  r \atop
{\scriptstyle i_{1},i_{2},\ldots\in M_{0} \atop \scriptstyle 
i_{1} < \ldots < i_{r} }}
(i_{1},\ldots,i_{r}) \in T_{r}.  \leqno(2)
$$
\end{theo}

The rest of this section is devoted to the proof of 2.1.  The key
will be the following

\begin{defi} Let  $M \subset \N$ be infinite, 
$i_{1},\ldots,i_{k} \in \N$, $ i_{1}< \ldots < i_{k} $. Moreover let
$(T_{r})_{r=1,2,\ldots}$  be as in 2.1.
\begaeq
\item $i_{1},\ldots,i_{k} \downarrow M$ abbreviates the
following fact:\\
If  $i_{k+1} < i_{k+2}< \ldots$  are points in $M$   with  $i_{k+1} >
i_{k}$,  then there is an  $r$  such that  $(i_{1},\ldots,i_{r})
\notin T_{r} $.  \\
The case  $k = 0$  will also be admissible, we will then write 
$\emptyset \downarrow M $.
\item $i_{1},\ldots,i_{k} \uparrow M$  stands for the
following:\\ 
Whenever  $N$  is an infinite subset of  $M $,  there are 
$i_{k+1} < i_{k+2} < \ldots$  in  $N$  with  $i_{k+1} > i_{k}$
such that $(i_{1},\ldots,i_{r}) \in T_{r}$  for every  $r $.  
\\
Again the 
definition is meant to contain the case  $k= 0 $.
\endaeq
\end{defi}
{\em Note.} Our  ``~$\downarrow$~''  and  ``~$\uparrow$~''  are closely
related with ``acceptance'' and  ``rejection'' in Diestel's book ([2],
p. 192). Our approach, however, is more direct since we only have
in mind a special version of a Ramsey theorem.

\bigskip
The following facts are immediate consequences of the definitions, they
are stated only for the sake of easy reference.

\bigskip\noindent
{\bf Observation 2.3}
\begaeq
\item (1) of the theorem just means  $\emptyset \uparrow \N $,  and
(2) follows as soon as one has found an  $M_{0}$  such that 
$i_{1},\ldots,i_{k} \uparrow M_{0}$  for arbitrary  $i_{1}< \ldots <
i_{k}$  in  $M_{0}$.
\item Let  $i_{1}< \ldots < i_{k}$  be given and  $M \subset
\N$  be infinite.  If  $i_{1},\ldots,i_{k} \downarrow M$ does not
hold, then there are  $i_{k+1} < i_{k+2} < \ldots$   (with  $i_{k+1}
> i_{k}$)  in  $M$  such that  $(i_{1},\ldots,i_{r}) \in T_{r}$ for
every  $r $.
\item If   $i_{1},\ldots,i_{k} \downarrow M$  and  $N \subset
M$  is infinite, then  $i_{1},\ldots,i_{k} \downarrow N $.  The same
holds if the ``~$\downarrow$~''  are replaced by ``~$\uparrow$''.
\endaeq
\bigskip

Here is the first step of our construction:

\setcounter{theo}{3}%
\begin{lemma}
There is an infinite  $\widetilde{M_{0}} \subset \N$ 
such that 
$$ i_{1},\ldots,i_{k} \downarrow \widetilde{M_{0}}
\quad{\rm\em or}\quad i_{1},\ldots,i_{k} \uparrow \widetilde{M_{0}}
\leqno(3)
$$
for each choice of  $i_{1} < \ldots < i_{k}$  in  $\widetilde{M_{0}}$ 
(including the case  $k=0$).
\end{lemma}
\proof It will be convenient to say that an infinite 
$\widetilde{M_{0}}$  satisfies  $(3)_{s}$  (where  $s \in \N_{0}$) if
the assertion (3) holds under the additional assumption that 
$\{i_{1},\ldots,i_{k}\}$  is contained in the set of the first  $s$ 
elements of   $\widetilde{M_{0}} $. We combine the following
observations:
\begpu
\item[--] We are looking for an infinite  $\widetilde{M_{0}}$  such
that  $\widetilde{M_{0}}$  satisfies (3)$_{s}$  for every  $s $.
\item[--] Suppose we are able to make the following induction work:
Given an infinite  $\widetilde{M}^{(s)}$  (which we write in
increasing order as  $\widetilde{M}^{(s)} =
\{i_{1},i_{2},\ldots,i_{s},\ldots\})$  such that  $(3)_{s}$  holds
for  $\widetilde{M}^{(s)} $,  there is an infinite subset  $N$  of 
$\{i_{s+1},i_{s+2},\ldots\}$  such that $\widetilde{M}^{(s+1)}:=
\{i_{1},\ldots,i_{s}\} \cup N$  satisfies   $(3)_{s+1} $.

This would suffice: We start our construction by setting 
$\widetilde{M}^{(0)}:= \N$  (note that  $\widetilde{M}^{(0)}$ 
satisfies  $(3)_{0}$  by  2.3(i)), use the induction to construct
the  $\widetilde{M}^{(0)} \supset \widetilde{M}^{(1)} \supset
\widetilde{M}^{(2)}\ldots$  and set  $\widetilde{M_{0}} =$ ``the
collection of the $s$'th elements of  $\widetilde{M}^{(s)}$,  $s \in
\N$.''  For fixed  $s$, $ \widetilde{M_{0}}$  has the same first  $s$ 
elements as  $\widetilde{M}^{(s)}$, so that in view of 2.3(iii)
 $(3)_{s}$ necessarily  holds for  $\widetilde{M_{0}} $.
\endpu
Therefore let's concentrate on the induction step.  Let  $s \geq 0$ 
and  $\widetilde{M}^{(s)}$  with  $(3)_{s}$  be given.  Denote by 
$\Delta_{1},\ldots,\Delta_{2^{s}}$  the  $2^{s}$  different subsets of 
$\{i_{1},\ldots,i_{s}\}$. 

We will construct infinite subsets  $\widetilde{N}^{[1]} \supset
\widetilde{N}^{[2]} \supset \ldots \supset \widetilde{N}^{[2^{s}]}$ 
of  $\{i_{s+2},\ldots\}$  such that either  $\Delta_{j},i_{s+1}
\uparrow \widetilde{N}^{[j]}$  or  $\Delta_{j},i_{s+1} \downarrow
\widetilde{N}^{[j]}$   for every  $j $.  In view of 2.3(iii) it is
then clear that  $\widetilde{M}^{(s+1)}:= \{i_{1},\ldots,i_{s+1}\}
\cup \widetilde{N}^{[2^{s}]}$  has  $(3)_{s+1} $.

First consider  $\Delta_{1}, i_{s+1} $.  Either we have  $\Delta_{1},
i_{s+1} \downarrow \widetilde{M}^{(s)}$  (in which case we put 
$\widetilde{N}^{[1]}:= \{i_{s+2},\ldots\}$)  or there is an infinite
subset  $\widetilde{N}^{[1]}$  of  $ \{i_{s+2},\ldots\}$   such that 
$\Delta_{1}, i_{s+1} \uparrow \widetilde{N}^{[1]}$  (see 2.3(ii)). 

Secondly, we investigate  $\Delta_{2}, i_{s+1} $.  Either 
$\Delta_{2},i_{s+1} \downarrow \widetilde{N}^{[1]}$  (we will put 
$\widetilde{N}^{[2]}:= \widetilde{N}^{[1]}$  in this case) or there is
an infinite subset  $\widetilde{N}^{[2]}$  of $\widetilde{N}^{[1]}$  such 
that   $\Delta_{2}, i_{s+1} \uparrow \widetilde{N}^{[2]} $.  It
should be clear how to construct the remaining  $\widetilde{N}^{[3]}
\supset \ldots \supset \widetilde{N}^{[2^{s}]} $.

In order to get an  $M_{0}$  with (2) from $\widetilde{M_{0}}$   we
need

\begin{lemma}
Let  $i_{1}< \ldots < i_{k}$  in 
$\widetilde{M_{0}}$  be given and suppose that  $i_{1},\ldots,i_{k}
\uparrow \widetilde{M_{0}} $.  Then there are only finitely many  $i
> i_{k}$  in  $\widetilde{M_{0}}$  such that  $i_{1},\ldots,i_{k}, i
\downarrow \widetilde{M_{0}} $.
\end{lemma}

\proof Suppose that this were not the case.  Put  $N = $  the
collection of these  $i $.  $N$  is infinite, and by 
$i_{1},\ldots,i_{k}  \uparrow \widetilde{M_{0}}$  there would be 
$i_{k+1} < i_{k+2} < \ldots $  in  $N$  (with  $i_{k+1}
> i_{k}$) such that  $(i_{1},\ldots,i_{r}) \in T_{r}$  for  every  $r
$.  Note that this would contradict  $i_{1},\ldots,i_{k+1} \downarrow
\widetilde{M_{0}} $.

\bigskip
Finally, we are ready for the

\bigskip\noindent
{\em Proof of theorem 2.1.} We have
already noted that (2) just means  $i_{1},\ldots,i_{k} \uparrow
M_{0}$   for  $i_{1} < \ldots < i_{k}$  in  $M_{0} $.  Similarly to
the proof of 2.4 we introduce  $(2)_{s}$  for  $s \geq 0$:  This is
(2) with the same additional assumption as in  $(3)_{s} $.

The construction parallels that of 2.4: We need an  $M_{0}$  with 
$(2)_{s}$  for all  $s$, we know that  $M^{(0)}:= \N$  satisfies 
$(2)_{0} $,  and, given an  $M^{(s)}$  with  $(2)_{s}$, 
we only have to construct an  $M^{(s+1)}$  with  
$(2)_{s+1}$.   In this
construction the first  $s$  elements of  $M^{(s)}$  and  $M^{(s+1)}$  
should be identical;  $M_{0} =$ ``the collection of the $s$'th elements
of  $M^{(s)}, s \in \N$'' then will have the desired
properties.

Here is the induction. Write  $M^{(s)} = \{i_{1},\ldots,i_{s},
i_{s+1},\ldots\}$  and put  $N:= \{i_{s+1},\ldots\} $.  Since 
$\Delta \uparrow N$  for every  $\Delta \subset
\{i_{1},\ldots,i_{s}\}$  by assumption we conclude from 2.5 that
there are only finitely many  $i$  in  $N$  such that  $\Delta,i
\downarrow N$  for any  $\Delta $.  Choose  $\widetilde{N} \subset N$ 
such that  $\widetilde{N}$  does not contain such  $i $.  Then 
$M^{(s+1)}:= \{i_{1},\ldots,i_{s}\} \cup \widetilde{N}$  satisfies 
$(2)_{s+1} $,  and this completes the proof.

\mysec{3}{A quantitative version of Rosenthal's $\ell^{1}$--theorem}%
To begin with, we restate a definition from [1].

\begin{defi}
Let  $(x_{n})$  be a bounded sequence in a Banach
space  $X$, and  $\varepsilon > 0. $  We say that  $(x_{n})$ admits 
$\varepsilon$--$\ell^{1}$--blocks  if for every infinite  $M\subset\N$ 
there are  $a_{1},\ldots,a_{r} \in \K$  with  $\sum |a_{r}| = 1$  and 
$i_{1} < \ldots < i_{r}$  in  $M$  such that  $\|\sum a_{\rho}
x_{i_{\rho}}\| \le \varepsilon $.  
\end{defi}

Clearly there will be no subsequence of   $(x_{n})$  equivalent to
the $\ell^{1}$--basis iff  $(x_{n})$  admits 
$\varepsilon$--$\ell^{1}$--blocks for arbitrarily small  $\varepsilon >
0 $.  Thus Rosenthal's theorem is the assertion that  $(x_{n})$  has
a weak Cauchy subsequence provided it admits 
$\varepsilon$--$\ell^{1}$--blocks for all  $\varepsilon $.  Here is our
quantitative version of this fact in the case of real spaces:

\begin{theo}
Let  $X$  be a real Banach space and 
$(x_{n})$  a bounded sequence.  Suppose that, for some  $\varepsilon
> 0 , (x_{n})$  admits small  $\varepsilon$--$\ell^{1}$--blocks.  Then
there is a subsequence  $(x_{n_{k}})$  of  $(x_{n})$  such that 
$(x_{n_{k}})$  is ``close to being a weak Cauchy sequence'' in the
following sense:
$$\limsup x'(x_{n_{k}}) - \liminf x'(x_{n_{k}}) \le
2 \varepsilon$$
for every  $x'$  with  $\|x'\| = 1 $.
\end{theo}
{\em Remark.} It is simple to derive the original theorem from 3.2. (If 
$(x_{n})$  and thus every subsequence has 
$\varepsilon$--$\ell^{1}$--blocks  for all  $\varepsilon $,  apply 3.2
successively with  $\varepsilon$ running through a sequence tending
to zero. The diagonal sequence which is obtained from this
construction will be a weak Cauchy sequence.) 
\bigskip
\proof Suppose the theorem were not true. We claim that 
without loss of generality
we may assume that there is a  $\delta > 0$  such that 
$$\varphi((x_{n_{k}})):= \sup_{\|x'\|=1} \Bigl( \limsup_{k}
x'(x_{n_{k}}) - \liminf_{k} x'(x_{n_{k}}) \Bigr) > 2\varepsilon +
\delta
\leqno(4)
$$
for all subsequences $(x_{n_{k}}) $.  In fact, if every subsequence
contained
another subsequence with a $\varphi$--value arbitrarily close to 
$2\varepsilon $,  an argument as in the preceding remark would even
provide one where  $\varphi((x_{n_{k}})) \le 2\varepsilon$  in
contrast to our assumption.

Fix a  $\tau > 0$  which will be specified later.  The essential tool
in order to get a contradiction will be the

\bigskip\noindent
{\bf Lemma} {\em After passing to a subsequence we may assume that 
$(x_{n})$  satisfies the following conditions:
\begaeq
\item If  $C$  and  $D$  are finite disjoint subsets of  $\N$ 
there are a  $\lambda_{0} \in \R$  and an  $x' \in X'$  with  $\|x'\|
=1$  such that  $x'(x_{n}) < \lambda_{0}$  for  $n \in C$  and 
$x'(x_{n}) > \lambda_{0} + 2\varepsilon + \delta$  for  $n \in D $.
\item There are  $i_{1} < \ldots < i_{r}$  in  $\N$, $
a_{1},\ldots,a_{r} \in \R$  with  
$$
\sum|a_{\rho}| = 1 ,\ |\sum a_{\rho}| \le
\tau,\ \|\sum a_{\rho} x_{i_{\rho}} \| \le \varepsilon .
$$
\endaeq
}

\bigskip\noindent
{\em Proof of the lemma.} 
(i) Define, for  $r \in \N , T_{r}$  to be the collection of all 
$(i_{1},\ldots,i_{r})$  (with  $i_{1} < \ldots < i_{r}$)  such that
there are a  $\lambda_{0} \in \R$  and a normalized  $x'$  such that 
$x'(x_{i_{\rho}}) < \lambda_{0}$  if  $\rho$  is 
even and  $>
\lambda_{0} + 2\varepsilon + \delta$  otherwise.  (4)
implies that (1) of 2.1 is valid. Thus there is an  $M_{0}$ 
for which all  $(i_{1},\ldots,i_{r})$  are in  $T_{r}$  for  $i_{1} <
\ldots < i_{r}$  in  $M_{0} $.  Let us assume that  $M_{0} = \N $.

Let  $C$  and  $D$  be finite disjoint subsets of  $2\N =
\{2,4,\ldots,\} $. We may select  $i_{1} < \ldots < i_{r}$  in  $\N$ 
such that  $C \subset \{i_{\rho}\mid\rho \; {\rm even}\}$  and  $D \subset
\{i_{\rho}\mid\rho \; {\rm odd}\} $.  Because of  $(i_{1},\ldots,i_{r})
\in T_{r}$  we have settled (i) provided  $C$  and  $D$  are in  $2\N
$,  and all what's left to do is to consider  $(x_{2n})$  instead of 
$(x_{n}) $.

(ii) By assumption we find  $i_{1} < \ldots < i_{r}$, $ a_{1},
\ldots, a_{r} \in \R $ such that $ \sum |a_{\rho}| = 1$ and $\|\sum
a_{\rho}x_{i_{\rho}}\| \le \varepsilon$  with arbitrarily large 
$i_{1} $.  Therefore we obtain  $i^{1}_{1} < \ldots < i^{1}_{r_{1}} <
i^{2}_{1} < \ldots < i^{2}_{r_{2}} < i^{3}_{1} < \ldots < i^{3}_{r_{3}}
< \ldots$  and associated  $a^{i}_{\rho} $.  The numbers 
$\eta_{j}:= \sum^{r_{j}}_{\rho=1} a^{j}_{\rho}$  all lie
in  $[-1,+1]$  so that we find  $j < k$  with  $|\eta_{j} - \eta_{k}| \le
2\tau $.  Let  $i_{1} < \ldots < i_{r}$  be the family  $i^{j}_{1} <
\ldots < i^{j}_{r_{j}} < i^{k}_{1} < \ldots < i^{k}_{r_{k}} $,  and
define the  $ a_{1},\ldots,a_{r}$  by  $\frac{1}{2}
a^{j}_{1},\ldots,\frac{1}{2} a^{j}_{r_{j}} , -\frac{1}{2}
a^{k}_{1},\ldots,-\frac{1}{2} a^{k}_{r_{k}} $.

\bigskip
We are now ready to derive a contradiction.  On the one hand, by (ii)
of the lemma, we find  $i_{1} < \ldots < i_{r}$  in  $\N $, $
a_{1},\ldots,a_{r} \in \R $, $\sum|a_{\rho}| = 1$, 
$|\sum a_{\rho}| \le \tau$ 
with  $\|\sum a_{\rho}x_{i_{\rho}}\| \le \varepsilon $.  On the other
hand  we may  apply (i) with  $C:= \{i_{\rho}\mid a_{\rho} < 0 \}$, $D:=
\{i_{\rho}\mid a_{\rho} > 0\} $.  We put  $\alpha:= -
\sum_{\rho\in C} a_{\rho}$, $ \beta := \sum_{\rho\in D}
a_{\rho} $,  and we note that  $|\alpha - \beta| \le \tau$, $\alpha +
\beta = 1$  so that  $|\beta - \frac{1}{2}| \le \tau $;  hence 
$$
\varepsilon \geq \left\|\sum a_{\rho}x_{i_{\rho}}\right\| 
\geq \sum a_{\rho}
x'(x_{i_{\rho}}) \geq - \lambda_{0}\alpha + (\lambda_{0} +
2\varepsilon + \delta)\beta \geq - |\lambda_{0}| \tau + \varepsilon +
\frac{\delta}{2} - \tau \delta .
$$ 
This expression can be made larger
than  $\varepsilon$  if $\tau$   has been chosen sufficiently small (note
that the numbers  $|\lambda_{0}|$  are bounded by  $\sup\|x_{n}\|) $, 
a contradiction which proves the theorem.

\bigskip\noindent
{\em Note.} Since  for the unit vector basis  $(x_{n})$  of real 
$\ell^{1}$  the assumption of the theorem holds with  $\varepsilon =
1$  and since for every subsequence  $(x_{n_{k}})$  one may find
$\|x'\| = 1$  with  $  \limsup x'(x_{n_{k}}) -
\liminf x'(x_{n_{k}}) =2$  
there can be no better constant than that given in
our theorem.

\bigskip

Let's now turn to the {\it complex case} which has to be treated in a
slightly different way. The quantitative version of Rosenthal's
$\ell^{1}$--theorem here reads as follows:

\begin{theo}
Let  $(x_{n})$  be a sequence in a
complex Banach space  $X$  such that, for some  $\varepsilon > 0 $, $
(x_{n})$  admits  $\varepsilon$--$\ell^{1}$--blocks.  Then there is a
subsequence  $(x_{n_{k}})$  such that, for  $x'$  with  $\|x'\| = 1 $, 
the diameter of the set of accumulation points of $(x'(x_{n_{k}}))_{k}$ 
is at most  $8\varepsilon/\sqrt{2} $.
\end{theo}
\proof The proof is similar to the preceding one. Again we know
that  -- if the theorem does not hold -- the numbers
$$\sup_{\|x'\|=1} \{\mbox{diameter of the accumulation points of }
(x'(x_{n_{k}}))\}$$
are greater than   $8\varepsilon/\sqrt{2} + \delta$ for a suitable 
$\delta > 0$  and all  $(x_{n_{k}}) $,  and again we fix a  $\tau > 0
$.

The key is this time the

\bigskip\noindent
{\bf Lemma} {\em  
Without loss of generality the sequence  $(x_{n})$ 
satisfies the following conditions.
\begaeq
\item Whenever  $C$ and $D$  are disjoint and finite subsets of  $\N$ 
there are  $z_{0}, w_{0} \in \C$  with  $|w_{0}| \geq
8\varepsilon/\sqrt{2} + \delta$  and an  $x'$ with  $\|x'\| = 1$ 
such that  $|x'(x_{n}) - z_{0}| \le \tau$  for  $n \in C$  and 
$|x'(x_{n}) - (z_{0} + w_{0})| \le \tau$ for $n \in D $.
\item We may assume that the  $a$'s   of 3.1 additionally
satisfy  $|\sum a_{\rho}| \le \tau$.
\endaeq
}

\bigskip
The proof is similar and is therefore omitted. 

Now let $S_{1},\ldots,S_{4}$ be the quadrants in the complex plane,
i.e.,
$S_{j}$  is the collection of those numbers whose arguments lie
between  $(j-1)\pi/2$  and  $j\pi/2 $.  We need the following obvious
facts:
\begpu
\item[--] If  $\sum |a_{\rho}| = 1$ there  is a  $j$  such that 
$\sum_{a_{\rho}\in S_{j}}|a_{\rho}| \geq 1/4$.
\item[--] If  $\sum_{a_{\rho}\in S_{j}}|a_{\rho}| \geq 1/4 $, 
then  $\sum_{a_{\rho}\in S_{j}} |a_{\rho}| \geq
\sqrt{2}/8 $.
\endpu
Now choose  $i_{1} < \ldots < i_{r}$, $a_{1}\ldots,a_{r}$  with  $
\sum|a_{r}| = 1$, $|\sum a_{\rho}| \le \tau$, $|\sum a_{\rho}
x_{i_{\rho}}| \allowbreak
\le \varepsilon $.  Write  $\{i_{1},\ldots,i_{r}\}$  as
the disjoint union of  $C$   and  $D $,   where  $D$ contains the 
$i_{\rho}$  with  $a_{\rho} \in S_{1}$  
(without loss of generality  we assume that 
$\sum_{a_{\rho}\in S_{1}}|a_{\rho}| \geq 1/4) $.

With $z_{0},w_{0},x'$  as in   (i)  of the lemma it follows that
\begin{eqnarray*}
\varepsilon ~\geq~ \left\|\sum a_{\rho}x_{i\rho}\right\| & \geq & 
\left|\sum
a_{\rho}x'(x_{i\rho}) \right|\\
& \geq & \left|\sum_{a_{\rho}\in S_{1}} a_{\rho}(z_{0}+w_{0})
+ \sum_{a_{\rho}\notin S_{1}} a_{\rho} z_{0}\right| - \tau
\sum|a_{\rho}|\\
& = & \left|\sum a_{\rho} z_{0} +\sum_{a_{\rho}\in S_{1}}
a_{\rho}w_{0}\right| -\tau\\
& \geq & |w_{0}| \; \left|\sum_{a_{\rho}\in S_{1}}a_{\rho}\right| -
\tau(1+|z_{0}|)\\
& \geq & |w_{0}| \dis\frac{\sqrt{2}}{8} - \tau(1+|z_{0}|)\\
& \geq & \varepsilon + \dis\frac{\sqrt{2}}{8} \delta -
\tau(1+|z_{0}|), 
\end{eqnarray*}
and this yields a contradiction for small  $\tau $.

\bigskip\noindent
{\em Note.} It seems to be difficult to decide whether  $8\eps/\sqrt{2}
$  can in general be replaced by a better estimate.

\mysec{4}{Blockings with uniformly bounded length}%
The following result has been published in [1].  Here we present a 
much simpler proof which uses theorem 2.1.
\begin{theo}
Let  $(x_{n})$  be a bounded sequence in a
real or complex Banach space  $X$  such that there is an  $r \in \N$ 
with the following property: Whenever  $M \subset \N$  is infinite
and  $\varepsilon > 0 $,  there are  $i_{1} < \ldots < i_{r}$  in  $M$ 
and  $a_{1},\ldots,a_{r} \in \K$  with  $\sum |a_{\rho}| = 1$ 
such that  $\|\sum a_{i} x_{i\rho}\| \le \varepsilon. $  Then 
$(x_{n})$  has a convergent subsequence.
\end{theo}

\proof  Fix  $\varepsilon > 0 $.  The first step is as in [1], we
refer the reader to this paper: One can choose the same  $a$'s  for
all  $M$.  Fix  $a_{1},\ldots,a_{r}$  and suppose that  $a_{r} \neq
0 $.  Define sets of  $\widetilde{r}$--tupels  $T_{\widetilde{r}}$ 
as follows.  $T_{\widetilde{r}}$  is the collection of {\it all} 
$\widetilde{r}$--tupels if  $\widetilde{r} \neq r $,  and the set of
those  $(i_{1},\ldots,i_{r})$  with  $i_{1} < \ldots < i_{r}$  and 
$\|\sum a_{\rho} x_{i\rho}\| \le \varepsilon $ if $\widetilde{r}=r$. 
Then 2.1 may be
applied, and we get  $(x_{n_{k}})$  such that  $\|\sum a_{\rho}
x_{n_{k_{\rho}}}\| \le \varepsilon$  for arbitrary
$k_{1} < \ldots < k_{r} $.  In particular   $x_{n_{r}} ,
x_{n_{r+1}} , x_{n_{r+2}}, \ldots$  lie in the ball with
center  $- \sum^{r-1}_{1} (a_{\rho}/a_{r}) x_{n_{\rho}}$  and
radius  $\varepsilon/|a_{r}| $.  Starting this construction with 
$\varepsilon = 1$  and applying it repeatedly to  $\varepsilon =
1/2$,
$\varepsilon = 1/3$, 
$\ldots$  provides a descending family of subsequences for 
which the diagonal sequence surely is convergent.

\bigskip\noindent
{\em Note.} In this situation a quantitative version is not to be
expected in general. Consider e.g. the unit vector basis  $(e_{n})$
in  $c_{0} $.  For  $\varepsilon > 0$  fixed one can find the 
$a_{1},\ldots,a_{r}$  with the same  $r$  for all  $M $,  but 
$\|e_{n} - e_{m}\| = 2$  for  $n \neq m $.

\mysec{5}{The Josefson--Nissenzweig theorem: preparations}%
\begin{lemma}
For  $\ell = 1,2,\ldots$  and 
$\varepsilon_{1},\ldots,\varepsilon_{\ell} \in \{0,1\}$  let 
$r_{\varepsilon_{1}\ldots \varepsilon_{\ell}}$  be a number such that
the family 
$((r_{\varepsilon_{1}\ldots\varepsilon_{\ell}})_{\varepsilon_{1}\ldots
\varepsilon_{\ell}})_{\ell=1,2,\ldots}$ 
satisfies  $|r_{\varepsilon_{1}\ldots \varepsilon_{\ell}}| \le
2^{-\ell}$  and  $r_{\varepsilon_{1}\ldots\varepsilon_{\ell}} =
r_{\varepsilon_{1}\ldots\varepsilon_{\ell}0} +
r_{\varepsilon_{1}\ldots\varepsilon_{\ell}1}$  for all  $\ell,
\varepsilon_{1}\ldots \varepsilon_{\ell} $.  Define  
$$\eta_{\ell}:=
\sum_{\varepsilon_{1}\ldots\varepsilon_{\ell}}
(r_{\varepsilon_{1}\ldots\varepsilon_{\ell}0} -
r_{\varepsilon_{1}\ldots\varepsilon_{\ell}1}) .
$$   
Then  $\sum
|\eta_{\ell}|^{2} \le 1$ so that in particular $\eta_{\ell}\to 0$.
\end{lemma}

\proof An elegant proof could be given using martingales: the  
$r_{\varepsilon_{1}\ldots\varepsilon_{\ell}}$  give rise to a bounded
martingale, the martingale convergence theorem guarantees the
existence of a limit  $f$  in  $L^{1} $,  and the  $\eta_{\ell}$ are the
integrals over  $f$  multiplied by suitable Rademacher functions.

However, a much simpler approach is possible.  Set
\begin{eqnarray*}
a_{\ell} & := & 2^{\ell}
\sum_{\varepsilon_{1}\ldots\varepsilon_{\ell}}
|r_{\varepsilon_{1}\ldots\varepsilon_{\ell}}|^{2}\\
b_{\ell} & := & 2^{\ell}
\sum_{\varepsilon_{1}\ldots\varepsilon_{\ell}}
|r_{\varepsilon_{1}\ldots\varepsilon_{\ell}0} -
r_{\varepsilon_{1}\ldots\varepsilon_{\ell}1}|^{2}\\
c_{\ell} &:= &
\sum_{\varepsilon_{1}\ldots\varepsilon_{\ell}}|r_{\varepsilon_{1}
 \ldots_{\ell}0} - r_{\varepsilon_{1}\ldots\varepsilon_{\ell}1}| .
\end{eqnarray*}
Then  $a_{\ell+1}- a_{\ell} = b_{\ell}$  (by the parallelogram law  
$|\alpha + \beta |^{2} + |\alpha-\beta|^{2} = 2(|\alpha|^{2} +
|\beta|^{2})$)  so that the   $a_{\ell}$ are increasing. Surely 
$a_{\ell} \le 1 $,  and we get  $\sum b_{\ell} \le 1 $. Finally
note that  $|\eta_{\ell}| \le c_{\ell}$  and that  $c^{2}_{\ell} \le
b_{\ell}$   since
$$(|\alpha_{1}| + \cdots + |\alpha_{k}|)^{2} \le k (|\alpha_{1}|^{2} +
\cdots + |\alpha_{k}|^{2}) \quad\mbox{for all families }
\alpha_{1},\ldots,\alpha_{k}.$$

Next we have to remind the reader of the definition of a Banach
limit  which can be found in nearly every textbook on functional
analysis. A {\em Banach limit\/} is an  $L \in
(\ell^{\infty})'$  such that   $L((1,1,\ldots)) = 1 = \|L\| $,  and 
$L((x_{1},x_{2},\ldots)) = L ((x_{2},x_{3},\ldots))$  for every 
$(x_{n}) $.  It is known that such  $L$  exist.  Some elementary
properties will be important for us:

\begin{lemma}
Define $\lambda_{\ell}$,
$\mu_{\varepsilon_{1}\ldots\varepsilon_{\ell}}\in \ell^{\infty}$  
for  $\ell = 1,2,\ldots$
and  $\varepsilon_{1},\ldots,\varepsilon_{\ell} \in \{0,1\}$  as
follows:
\begin{eqnarray*}
\lambda_{1} & = & (1, -1, 1, -1, \ldots)\\
\lambda_{2} & = & (1, 1, -1, -1, 1, 1, \ldots)\\
\lambda_{3} & = & (1,1,1,1,-1,-1,-1,-1,\ldots)\\
&\vdots &\\
\mu_{0} & = & (1,0,1,0,\ldots),\quad \mu_{1}~ =~ (0,1,0,1,\ldots)\\
\mu_{00}  & = & (1,0,0,0,1,0,0,0,1,\ldots),\quad \mu_{01} ~=~
(0,0,1,0,0,0,1,0,\ldots)\\
\mu_{10} & = & (0,1,0,0, 0,1, 0,0,0,\ldots),\quad \mu_{11} ~=~
(0,0,0,1,0,0,0,1,\ldots);
\end{eqnarray*}
in general: $\mu_{\varepsilon_{1}\ldots\varepsilon_{\ell}}$  is  $1$ 
at positions of the form  $k \cdot 2^{\ell} + 1 + \varepsilon_{1}2^{0}
+ \cdots + \varepsilon_{\ell} 2^{\ell-1}$ 
$(k=0,1,\ldots)$  and  $0$ 
otherwise, and  
$$\lambda_{\ell} =
\sum_{\varepsilon_{1},\ldots,\varepsilon_{\ell-1}}
(\mu_{\varepsilon_{1},\ldots,\varepsilon_{\ell-1},0} -
\mu_{\varepsilon_{1},\ldots,\varepsilon_{\ell-1},1} ).
$$
Then  $|L(\mu_{\varepsilon_{1}\ldots\varepsilon_{\ell}} x)| \le 2^{-\ell}$ 
for every  $x \in \ell^{\infty}$  with  $\|x\| \le 1 $.
\end{lemma}
\proof  Let   $T: \ell^{\infty} \to \ell^{\infty}$  be the
shift operator  $(y_{1},y_{2},\ldots)\mapsto (0,y_{1},y_{2},\ldots)$ and 
$x_{0}$  the pointwise product of 
$\mu_{\varepsilon_{1}\ldots\varepsilon_{\ell}}$  with  $x $.  Then 
$L(T\widetilde{x}) = L(\widetilde{x})$  for every  $\widetilde{x} $,
and  $\|x_{0} + T x_{0} + T^{2}x_{0} + \cdots + T^{2^{\ell}-1}x_{0}\|
\le 1 $.   Hence  $2^{\ell}|L(x_{0})|\le 1 $.

\mysec{6}{The Josefson--Nissenzweig theorem}%
\begin{theo}[{[}4{]}, {[}5{]}]  
A Banach space  $X$  is either
finite--dimensional or there exist normalized  $x'_{1},
x'_{2},\ldots$  such that  $x'_{n}(x) \to 0$  for every  $x $.
\end{theo}
\proof  \\
{\sc Case 1:} $\ell^{1}$  is not contained in  $X'$. Suppose that 
every weak*--convergent sequence is already
norm convergent; we will show that  $X'$  is
finite-dimensio\-nal.

Let  $(x'_{n})$  be a bounded sequence.  By Rosenthal's theorem and
since $\ell^{1}$  does not embed into  $X'$  we find a subsequence
which is weakly Cauchy and thus weak*--convergent.  By our assumption
it is convergent, and thus  $X'$  is finite--dimensional.

\bigskip\noindent
{\sc Case 2:}  $\ell^{1}$  embeds into  $X'$, i.e., there are 
$x'_{n}$  in  $X'$  and  $A,B > 0$  such that  
$$A \sum^{r}_{1}
|t_{i}| \geq \left\|\sum^{r}_{1} t_{i} x'_{i}\right\| \geq B
\sum^{r}_{1}|t_{i}|
$$  
for arbitrary  $r$  and 
$t_{1},\ldots,t_{r} \in \K $.

In order to continue we remind the reader of the following notion: A
sequence  $(y'_{n})$  is said to be obtained from the  $(x'_{n})$  by
{\em blocking\/}  if there are disjoint finite sets 
$A_{1},A_{2},\ldots$  in  $\N$  with  $A_{1} \le A_{2} \le \ldots$  and
numbers  $(a_{k})$  with  $\sum_{k\in A_{n}}|a_{k}| = 1$  for
every  $n$  such that  $y'_{n} = \sum_{k \in A_{n}} a_{k}
x'_{k} $.
Note in particular that all subsequences arise in this way.

\medskip
\noindent
{\sc Case 2.1.}  It is possible to get  $(y'_{n})$ by blocking
$(x_{n}) $ such that  $y'_{n}
\rightarrow 0$  with respect to the weak*--topology.  Then we are done
since  $\|y'_{n}\| \geq B$ so that the  $y'_{n}/\|y'_{n}\|$  have the
desired properties.

\medskip
\noindent
{\sc Case 2.2.} For no blocking  $(y'_{n})$  we have  $y'_{n}
\rightarrow 0$  (w.r.t. the weak*--topology).  
In order to measure the property 
of being a weak*--null sequence we introduce the number 
$$\varphi((y'_{n})) := \sup_{\|x\|=1}  \limsup |y'_{n}(x)|$$
for the  $(y'_{n})$  constructed as before.  In the case under
consideration we know that always  $\varphi((y'_{n})) > 0 $.  We
claim that even more is true.

\medskip
\noindent
{\sc Claim 1.} There are a  $\delta > 0$  and a block sequence 
$(y'_{n})$  such that  $\varphi((z'_{n})) = \delta$  for every 
$(z'_{n})$  which is obtained from  $(y'_{n})$  by blocking.

\medskip
\noindent
{\em Proof of claim 1.}  Let  $\delta_{0} \geq 0$  be the  infimum  of the
numbers  $\varphi((y'_{n})) $, where the  infimum  runs over all block
sequences  $(y'_{n}) $.  Choose  $(y^{[1]}_{n}), $  a block sequence
of  $(x'_{n}) $,  such that  $\varphi((y^{[1]}_{n})) \le \delta_{0} +
1/2^{0} $. Let  $\delta_{1}$  be the infimum  of the 
$\varphi((y'_{n})) $,   where this time only blockings of 
$(y^{[1]}_{n})$  are under consideration.  Since a block sequence of
a block sequence is
a block sequence we have  $\delta_{0} \le \delta_{1} $.  Choose 
$(y^{[2]}_{n}) $, a block  sequence of  $(y^{[1]}_{n}) $,  with 
$\varphi((y^{[2]}_{n})) \le \delta_{1} + 1/2^{1} $.  In this way  we get
successively  $(y^{[1]}_{n}) , (y^{[2]}_{n}), \ldots$,  and 
$\delta_{1} \le \delta_{2} \le \ldots $,  where   $(y^{[k+1]}_{n})$ 
is obtained from  $(y^{(k)}_{n})$ by blocking and where  $\delta_{k}
\le \varphi((y'_{n})) \le \delta_{k} + 1/2^{k}$  for all block
sequences  $(y'_{n})$ of   $(y^{[k+1]}_{n}) $.  Our candidate is the
diagonal sequence  $(y'_{n})$  containing the $n$'th element of the
$(y^{[n]}_{n})$  for every  $n $.  
For every $k$, $(y_{n}')$ is -- after possibly finitely many
exceptions -- a block sequence of $(y_{n}^{[k]}$; therefore the 
$\varphi$--value lies between  $\delta_{k+1}$  and  $\delta_{k} +
1/2^{k} $.  It follows that  $\delta := \sup \delta_{k}$  has the
claimed properties; note that we also know that  $\delta > 0$  since 
$\delta = \varphi((y'_{n})) $.
                        
\medskip
\noindent
{\sc Claim 2.}  Fix  $(y'_{n})$  and  $\delta > 0$  as in claim 1. 
Further let  $\tau > 0$  be arbitrary.  There is a subsequence which
we will denote by  $(z'_{n})$  with the following property: It is
possible to find normalized  $x_{1},x_{2},\ldots$  in  $X$  such
that:
\begin{eqnarray*}
 |z'_{n}(x_{1}) - \delta| &\le& \tau   \mbox{ for all } 
n;\\
|z'_{n}(x_{2}) - \delta| &\le& \tau   \mbox{ for } 
n=3,5,\ldots\\
|z'_{n}(x_{2}) + \delta| &\le& \tau   \mbox{ for } 
n=4,6,\ldots\\
|z'_{n}(x_{3}) - \delta| &\le& \tau   \mbox{ for } 
n=5,6,9,10,\ldots\\
|z'_{n}(x_{3}) + \delta| &\le& \tau   \mbox{ for } 
n=7,8,11,12,\ldots\\
&\vdots &
\end{eqnarray*}
(In general  $z'_{n}(x_{k})$  is $\tau$--close to  $\delta$  on
segments of length  $2^{k-2}$ beginning at  $2^{k-1}+1 , 2 \cdot
2^{k-1}+1,3 \cdot 2^{k-1}+1, \ldots$  and $\tau$--close to  $- \delta$ 
at segments of the same length, beginning at  $2^{k-1} + 2^{k-2}+1$, $2
\cdot 2^{k-1} + 2^{k-2}+1$, $3 \cdot 2^{k-1} + 2^{k-1}+1$;  
thus, if we regard the 
$x_{k}$  as functions on the set  $\{z'_{n}\mid n \in \N\}$ they behave
like the Rademacher functions, at least for large $n$.)

\medskip
\noindent
{\em Proof of claim 2.} By assumption we know that 
$\varphi((y'_{n})) = \delta $,  and this makes it easy to find  $x_{1}$ 
with  $\|x_{1}\| = 1$  and a subsequence  $(w^{[1]}_{n})$  of 
$(y'_{n})$  with  $|w^{[1]}_{n}(x_{1}) - \delta| \le \tau$  for every 
$n $.  Our final  $(z'_{n})$  will be a subsequence of $(w^{[1]}_{n})$
so that we will have no problems with  $x_{1} $.  Put  $z'_{1} :=
w^{[1]}_{1} $, $z'_{2}:= w^{[2]}_{2} $.
Now consider  
$$(u_{n}) := \left(\dis\frac{w^{[1]}_{3}-w^{[1]}_{4}}{2}
, \frac{w^{[1]}_{5}-w^{[1]}_{6}}{2}, \ldots\right) .$$
This is a block sequence of  $(y'_{n})$ so that  $\varphi((u_{n})) =
\delta $.  Hence we find a normalized  $x_{2}$  such that for
infinitely many  $n$,  say  $n \in N $,  we have  $|u_{n}(x_{2}) -
\delta| \le \tau' $;  here  $\tau'$  denotes any positive number such
that  
$$
|\alpha|,|\beta| \le \delta + \tau',\ 
|\frac12 (\alpha+\beta) -  \delta | \le \tau'\quad \Rightarrow \quad
|\alpha-\delta|\le\tau,\   |\beta-\delta| \le \tau .
$$  
Since  
$\limsup |w^{[1]}_{n}(x_{2})| \le \delta$  we may also assume that 
for $n\in N$  and  $u_{n}= (w^{[1]}_{2n-1} - w^{[1]}_{2n})/2 $,  the 
$w$'s  satisfy  $|w^{[1]}_{2n-1}(x_{2})|$, $ |w^{[1]}_{2n}(x_{2})| \le
\delta + \tau' $. Thus  $|w^{[1]}_{2n-1}(x_{2}) - \delta| \le \tau$ 
and  $|w^{[1]}_{2n}(x_{2}) + \delta| \le \tau $.

Let  $(w^{[2]}_{n})$  be the sequence of the  $w^{[1]}_{2n-1} ,
w^{[1]}_{2n}$  with  $n \in N $.  Set  $z'_{3}:= w^{[2]}_{2} $, $
z'_{4}:= w^{[2]}_{3} $. 

We consider now  
$$ \left(\dis\frac{w^{[2]}_{1}+w^{[2]}_{2} -
w^{[2]}_{3}-w^{[2]}_{4}}{4} , \frac{w^{[2]}_{5} + w^{[2]}_{6} -
w^{[2]}_{7} - w^{[2]}_{8}}{4} , \ldots\right) =: (u_{n}) .
$$ 
 Again we
find an infinite  $N$  and an  $x_{3}$  such that 
$|u_{n}(x_{3})-\delta| \le \tau' $, where this time  $\tau' > 0$  is
such that  
$$
|\alpha_{1}|, |\alpha_{2}|, |\alpha_{3}|, |\alpha_{4}|
\le \delta + \tau' 
\mbox{ and }
 \left|\dis\frac{\alpha_{1} +
\alpha_{2}+\alpha_{3}+\alpha_{4}}{4} - \delta\right| \le \tau'
$$
%imply 
$$
\Rightarrow\quad
|\alpha_{j} - \delta| \le \tau \mbox{ for }j = 1,2,3,4 .
$$  
Also we can
assume that the  $|w^{[2]}_{4n-j}(x_{3})| \le \delta + \tau'$  for 
$n \in N$, $ j = 0,1,2,3 $. Let   $(w^{[3]})$  consist of the  
$w^{[2]}_{4n-3},\ldots,w^{[2]}_{4n}$  with  $n \in N
$,  and define  $z'_{5},\ldots,z'_{8}$ to be 
$w^{[3]}_{1},\ldots,w^{[3]}_{4}$  respectively.

It should be clear how this construction (which is similar to that in
[3]) has to be continued and
that  $(z'_{n})$  has the claimed properties.

Now it is fairly easy to conclude the proof of the
Josefson--Nissenzweig theorem.  With  $(z'_{n})$  as in the second
claim we denote by  $T: X \to \ell^{\infty}$  the operator  $x
\mapsto (z'_{1}(x), z'_{2}(x),\ldots) $,  and we define  $w'_{n}$  to
be the functional  $x \mapsto L (\lambda_{n}Tx)$  (notation as in
section 5) where $L$ is a fixed Banach limit.  
Every  $w'_{n}$  has a norm not smaller than $\delta -
\tau$  since  $\lambda_{n} Tx_{n}$ is by construction a sequence which  -- up
to finitely many exceptions -- is $\tau$--close to the sequence 
$(\delta, \delta, \delta,\ldots) $,  and the Banach limit property of 
$L$  implies that  $L(\lambda_{n} Tx_{n})$  is -- up to $\tau$ -- 
$L(\delta,\delta,\ldots) = \delta $.

The  $w'_{n}$  also tend to zero w.r.t. the weak* topology since by 5.2 
for every  $x \in \ell^{\infty}$ with  $\|x\| = 1$  the numbers 
$r_{\varepsilon_{1}\ldots\varepsilon_{\ell}} :=
L(\mu_{\varepsilon_{1}\ldots\varepsilon_{\ell}} x)$  satisfy the
hypothesis of 5.1 and since the 
$\eta_{\ell} $ of 5.1 are just the 
$L(\lambda_{\ell}\cdot x)$ in this case.

\bigskip\noindent
{\em Note.} In fact we have shown more than required:  When case
2.2 leads to a Josefson--Nissenzweig sequence then it is not only
weak* null but in fact a weak*--$\ell^{2}$--sequence (i.e. 
$\sum|w'_{n}(x)|^{2} < \infty$  for every  $x \in X) $.

\begin{center}
{\sc References}
\end{center}
\begingroup
\renewcommand{\labelenumi}{\mbox{\rm[\arabic{enumi}]}}
\begin{enumerate}
\item
{Behrends, E.}: On Rosenthal's $\ell^{1}$--theorem\\ 
to appear (Archiv d. Math., 1994).
\item
{Diestel, J.}: Sequences and Series in Banach Spaces\\
Springer--Verlag, Berlin--Heidelberg--New York, 1984.
\item
{Hagler, J.} and Johnson, W.B.: On Banach spaces whose dual
balls are not weak* sequentially compact\\
Israel J. Math. 28, 1977,
325--330.
\item
{Josefson, B.}: Weak sequential convergence in the dual of a
Banach space does not imply norm convergence\\
Ark. Mat. 13, 1975, 79--89.
\item
{Nissenzweig, A.}: $w*$ sequential convergence\\
Israel J. Math. 22, 1975, 266--277.
\item
{Rosenthal, H.P.}: A characterization of Banach spaces
containing $\ell^{1}$\\
Proc. Nat. Acad. Sci. (USA) 71, 1974, 2411--2413.
\end{enumerate}
\endgroup

\bigskip
\bigskip
\bigskip\noindent
I. Mathematisches Institut\\
FU Berlin\\
Arnimallee 2--6\\
D--14195 Berlin\\
e-mail: behrends@math.fu-berlin.de

\end{document}